\documentclass{amsart}

\usepackage{ amsmath,amssymb,amsbsy,amsfonts, latexsym,amsopn,amstext, amsxtra,euscript,amscd, hyperref}
\usepackage{graphicx}
\usepackage{cite}

\usepackage{amsrefs}

\usepackage{url}

\newcommand\myurl[1]{\url{#1}}

\BibSpec{webpage}{%
  +{}{\PrintAuthors} {author}
  +{,}{ \textit} {title}
  +{}{ \parenthesize} {date}
  +{,}{ \myurl} {myurl}
  +{,}{ } {note}
  +{.}{ } {transition}
}

\usepackage{xy} 
\input xy \xyoption{all}

\newtheorem{thm}{Theorem}
\newtheorem{prop}{Proposition}
\newtheorem{lemma}{Lemma}
\newtheorem{rem}{Remark}

\newtheorem{exa}{Example}

\newtheorem{defn}{Definition}

\def\Z{\mathbb Z}
\def\Q{\mathbb Q}
\def\C{\mathbb C}
\def\P{\mathbb P^1 (k)}
\def\Q{{\mathbb Q}}
\def\Z{{\mathbb Z}}
\def\C{{\mathbb C}}
\def\L{\mathcal L}
\def\H{\mathcal H}
\def\M{\mathcal M}
\def\X{\mathcal X}
\def\S{\mathcal S}

\def\p{\mathfrak p}

\def\p{\mathfrak p}
\def\s{\mathfrak s}

\def\Jac{\mbox{Jac }}

\def\D{\Delta}

\def\emb{\hookrightarrow}

\def\a{\alpha}
\def\e{\varepsilon}

\def\t{\tau}
\def\<{\langle}
\def\>{\rangle}
\def\zz{\zeta}
\def\z{\omega}

\def\iso{{\, \cong\, }}

\def\bG{\bar G}

\def\bG{{\bar G}}

\def\sem{{\rtimes}}

\def\Z{\mathbb Z}

\def\Q{\mathbb Q}

\def\C{\mathbb C}

\def\C{\mathcal C}
\def\H{\mathcal H}
\def\M{\mathcal M}

\def\a{\alpha}

\def\p{\mathfrak p}

\def\Z{\mathbb Z}
\def\Q{\mathbb Q}
\def\C{\mathbb C}
\def\L{\mathcal L}
\def\H{\mathcal H}
\def\M{\mathcal M}

\def\s{\sigma}
\def\a{\alpha}
\def\p{\mathfrak p}

\def\<{\langle}
\def\>{\rangle}
\def\X{\mathcal X}
\def\emb{\hookrightarrow}

\def\D{\Delta}

\def\L{\mathcal L}

\def\S{\mathcal S}

\def\Y{\mathcal Y}

\def\P{\mathbb P^1 (k)}
\def\Q{{\mathbb Q}}
\def\Z{{\mathbb Z}}
\def\C{{\mathbb C}}
\def\M{{\mathcal M}}
\def\H{\mathcal H}
\def\L{\mathcal L}

\def\J{\mbox{Jac }}

\def\Aut{\mbox{Aut }}
\def\bAut{\overline {\mathrm{Aut}}}

\def\Jac{\mbox{Jac }}

\def\embd{\hookrightarrow}

\def\D{\Delta}
\def\e{\varepsilon}

\def\bG{\bar G}
\def\t{\tau}

\def\a{{\alpha }}

\def\<{\langle}
\def\>{\rangle}

\def\sem{{\rtimes}}

\def\s{\mathfrak s}

\def\p{\mathfrak p}

\def\emb{\hookrightarrow }

\def\zz{\zeta}

\def\z{\omega}

\hypersetup{
  colorlinks   = true, %Colours links instead of ugly boxes
  urlcolor     = blue, %Colour for external hyperlinks
  linkcolor    = blue, %Colour of internal links
  citecolor   = red %Colour of citations
}

\numberwithin{equation}{section}

\begin{document}

\title{Genus 3 hyperelliptic curves with $(2, 4, 4)$-- split Jacobians}

\author{T. Shaska}
\address{Department of Mathematics, Oakland  University, Rochester, MI, 48309}
%\curraddr{Department of Mathematics and Statistics,Case Western Reserve University, Cleveland, Ohio 43403}
\email{shaska@oakland.edu}

%    General info
\subjclass[2000]{Primary 20F70, 14H10; Secondary 14Q05, 14H37}

%\date{January 1, 2001 and, in revised form, June 22, 2001.}


\keywords{invariants, binary forms, genus 3, algebraic curves}

\begin{abstract}
We study degree 2 and 4 elliptic subcovers of  hyperelliptic curves of genus 3 defined over $\C$.   The family   %$\mathcal T_3$ 
of genus 3 hyperelliptic curves which have a degree 2 cover to an elliptic curve $E$ and degree 4 covers to  elliptic curves $E_1$ and $E_2$ is a 2-dimensional subvariety of the hyperelliptic moduli $\H_3$. We determine this subvariety explicitly. For any given moduli point $\p \in \H_3$ we determine explicitly if the corresponding genus 3 curve $\X$ belongs or not to such family. % $\mathcal T_3$. 
When it does, we can determine elliptic subcovers $E$, $E_1$, and $E_2$  in terms of the absolute invariants $t_1, \dots , t_6$ as in \cite{hyp_mod_3}.
This variety  provides a new family of hyperelliptic curves of genus 3 for which the Jacobians completely split. 
The sublocus  of such family %$\mathcal T_3$ 
when  $E_1 \iso E_2$ is a 1-dimensional variety which we determine explicitly.  We can also  determine $\X$ and $E$ starting form the $j$-invariant of $E_1$. 
\end{abstract}

\maketitle

\section{Introduction}

Let $\M_g$ denote the moduli space of genus $g\geq 2$ algebraic curves defined over an algebraically closed field $k$ and $\H_g$ the hyperelliptic submoduli in $\M_g$.  The sublocus of genus g hyperelliptic curves with an elliptic involution is a $g$-dimensional subvariety of $\H_g$.  For $g=2$ this space is denoted by $\L_2$ and studied in Shaska/V\"olklein \cite{sh-v} and   for $g=3$ is denoted by $\S_2$ and is computed and discussed in detail   in \cite{b-th}.  In both cases, a birational parametrization of these spaces is found via \textit{dihedral invariants} which are introduced by this author in \cite{issac, g-sh-s, sh_03} and generalized for any genus $g \geq 2$ in \cite{g-sh}.  We denote the parameters for $\L_2$ by $u, v$ and for $\S_2$ by $\s_2, \s_3, \s_4$ as in respective papers.   Hence, for the case $g=3$ there is a birational map 
\begin{equation}
\begin{split}
 \S_2 & \longrightarrow \H_3 \\
\phi : (\s_2, \s_3, \s_4) & \longrightarrow  (t_1, \dots , t_6), 
\end{split}
\end{equation}
where $t_1, \dots , t_6$ are the absolute invariants as defined in \cite{hyp_mod_3} satisfying the equation of the genus 3 hyperelliptic moduli; see \cite{hyp_mod_3} for details. 

The dihedral invariants $\s_2, \s_3, \s_4$ provide a birational parametrization of the locus $\S_2$.  Hence, a generic curve in $\S_2$ is uniquely determined by the corresponding triple $(\s_2, \s_3, \s_4)$.  

Let $\X$ be a curve in the locus $\S_2$.  Then there is a degree 2 map $f_1:   \X \to E$ for some elliptic curve $E$.  Thus, the Jacobian of $\X$ splits as $\Jac (\X) \iso E \times A$, where $A$ is a genus 2 Jacobian.  Hence, there is a map $f_2 : \X \to C$ for some genus 2 curve $C$.  The equations of $\X$, $E$, and $C$ are given in Thm.~\ref{thm_1}.  For any fixed curve $\X \in \S_2$, the subcovers  $E$ and $C$ are uniquely determined in terms of the invariants $\s_2, \s_3, \s_4$. 

In section three we give the splitting of the Jacobians for all genus 3 algebraic curves, when this splitting is induced by automorphisms.  The case for hyperelliptic curves was done in \cite{P}, we complete the non-hyperelliptic cases.  The proof requires the Poincare duality and some basic group theory. 

In this paper, we are mostly interested in the case when the Jacobian of the  genus two curve $C$ also splits.  The Jacobian of $C$ can split as an $(n, n)$-structure; see \cite{sh_01}.  The loci of such genus 2 curves with    $(3, 3)$-split  or $(5, 5)$-split have been studied respectively in in \cite{sh_01} and \cite{deg5}.  We focus on the case when the Jacobian of $C$ is $(2, 2)$-split, which corresponds to the case when the Klein 4-group $V_4 \emb \Aut (C)$.   Hence, $\Jac \X$ splits completely as a product of three elliptic curves. We say that $\Jac \X $ is $(2, 4, 4)$-split. 

Let the locus of genus 3 hyperelliptic curves whose Jacobian is $(2, 4, 4)$-split  be denoted by $\mathcal T$.    Then, there is a rational map $\psi: \mathcal T \to \L_2$ such that 
\begin{equation} \psi : (\s_2, \s_3, \s_4 ) \longrightarrow  (u, v) \end{equation}
which has degree 70 and can be explicitly computed, even though the rational expressions of $u$ and $v$ in terms of $\s_2, \s_3, \s_4$ are quite large. 

There are three components of $\mathcal T$ which we denote them by $\mathcal T_i$, $i=1, 2, 3$.  Two of these components are well known and the correspond to the cases when $V_4$ is embedded in the reduced  automorphism group of $\X$.  These cases correspond to the singular locus of $\S_2$ and are precisely the locus $\det \left( \Jac (\phi) \right)=0$, see \cite{sh_03}.  This happens for all genus $g\geq 2$ as noted by Geyer \cite{Ge}, Shaska/V\"olklein \cite{sh-v}, and shown in \cite{sh_03}.  The third component $\mathcal T_3$  is more interesting to us.  It doesn't seem to have any group theoretic reason for this component to be there in the first place.  We find the equation of this component it terms of the $\s_2, \s_3, \s_4$ invariants.  It is an equation 
\begin{equation} F_1 (\s_2, \s_3, \s_4 ) =0 \end{equation}
as in Eq.~\eqref{eq_T}. In this locus, the elliptic subfields of the genus two field $k(C)$ can be determined explicitly.  

\def\SV{\mathcal V}

The main goal of this paper is to determine explicitly the family $\mathcal T_3$ of genus 3 curves and relations among its elliptic subcovers.
We have the maps 
\begin{equation}
\begin{aligned}
  \mathcal T_3  &  \overset{\psi}{\longrightarrow} \L_2    \overset{\psi_0}{\longrightarrow}  k^2  \\
  (\s_2, \s_3, \s_4 )  &\to  (u, v)  \to (j_1, j_2) \\
\end{aligned}
\end{equation}
where $\s_2, \s_3, \s_4$ satisfy  Eq.~\eqref{eq_T} and $u, v$ are given explicitly by Eq.~\eqref{i1_i2_i3} and Thm.~ (3) in \cite{sh-v}. The degree $\deg \, \psi_0 =2$ and $\deg \, \psi =70 $.

Since $\mathcal T_3$ is a subvariety of $\H_3$ it would be desirable to express its equation in terms of a coordinate in $\H_3$.  One can use the absolute invariants of the genus 3 hyperelliptic curves $t_1, \dots , t_6$ as defined in \cite{hyp_mod_3} and the expressions of $\s_2, \s_3, \s_4$ in terms of these invariants as computed in \cite{b-th}.

Further, we focus our attention to the sublocus $\SV$ of $\L_2$ such that the genus 2 field $k(C)$ has isomorphic elliptic subfields.  Such locus was discovered by Shaska/V\"olklein in \cite{sh-v} and it is somewhat surprising.  It does not rise from a family of genus two curves with a fixed automorphism group as other families, see \cite{sh-v} for details. Using this sublocus of $\M_2$ we discover a rather unusual embedding $\M_1 \emb \M_2$ as noted in \cite{sh-v}.  Let $\mathfrak T \subset \mathcal T_3 \subset \H_3$  be the subvariety of $\mathcal T_3$ obtained by adding the condition $j_1=j_2$.   Then, $\mathfrak T$ is a 1-dimensional variety  defined by equation 
\begin{equation}
\left\{ 
\begin{split}
F_1 (\s_2, \s_3, \s_4) = 0 \\
F_2 (\s_2, \s_3, \s_4) = 0 \\
\end{split}
\right.
\end{equation}
where $F_2 $ is the discriminant of the quadratic polynomial roots of which are $j$-invariants $j_1$ and $j_2$; cf. Lemma~3.  
%The variety $\mathfrak T$ is a rational curve ????????.  
Hence, we have the maps
\begin{equation}
\begin{split}
k     \to     \mathfrak T          \emb      \SV    \emb   k          \\
  t   \to    (\s_2, \s_3, \s_4)   \to      (u, v)   \to       j_1       \\
\end{split}
\end{equation}

Next we study whether the above maps are invertible.  That would provide birational parametrizations for varieties $\SV$ and $\mathfrak T$.  The variety $\SV$ is known to have a birational parametrization from Thm.~ (3) in \cite{sh-v}.  The map can be inverted as 
\begin{equation}
\begin{aligned}
%\M_1 & \emb  & \SV                 & \\
 j   & \to   & (u, v)= \left(   9 - \frac j {256} ,   9 \left( 6 - \frac j {256}\right)  \right),  & \\
\end{aligned}
\end{equation}
see \cite{sh-v} for details. The main computational task of this paper is to find a birational parametrization of $\mathfrak T$.    
%Hence, we can invert the maps  as follows
%\[ 
%\begin{aligned}
%\M_1 & \emb  & \SV   & \emb  & \mathcal T               & \\
%  j  & \to   & (u, v) & \to   & (\s_2, \s_3, \s_4)  & \\
%\end{aligned}
%\]
%

Given $(u, v) \in \SV$ there is a unique (up to isomorphism) genus 2 curve $C$ corresponding to this point in $\SV$.  From Lemma~\ref{lem_gen_2}, every genus 2 curve can be written as in Eq.~\eqref{g3-g2}.  Hence, there exists a triple $(\s_2, \s_3, \s_4)$ corresponding to $(u, v)$.   

%Expressing $\s_2, \s_3, \s_4$ in terms of $j$ we find a rational parametrization of the variety $\mathcal T$. 

If $j \in \Q$, then it is a well-known fact that the corresponding elliptic curve $E_j$ can be chosen defined over $\Q$.  From the above expressions we see that $u, v \in \Q$.  Then, from   \cite{sh_02} the corresponding genus two curve $C$ has also minimal field of definition $\Q$.  The same holds for $\s_2, \s_3, \s_4$ and the genus 3 corresponding curve $\X$. 

%Hence, we have the above maps giving
%
%\[ E/\Q \to C/\Q \to \X/ \Q \]

%The rest of the paper compares the elliptic subfield $E$ with the two isomorphic subfields $E_1$ and $E_2$.  Each one of such subfields is determined uniquely (up to isomorphism).  We find necessary and sufficient conditions when $E$ is isomorphic with $E_1$ and when they are 3-isogenous and 5-isogenous.  

%The family $\mathcal T$ could not be discovered without using computational techniques, since they do not have a fixed automorphism group and would have been missed by going through the list of groups as  in \cite{b-th} or other papers \cite{g-sh-s, sh_03, issac}.  

%\bigskip
%\noindent\textbf{Notation}   Throughout this paper by a curve we mean an irreducible algebraic curve defined over an algebraically closed field of characteristic zero. 

%**********************************************************************
%\newpage
\section{Hyperelliptic curves with extra involutions}

In this section, we will go over some basic properties of genus three hyperelliptic fields which have an elliptic involution.  For details on genus 3 curves the reader can check further \cite{dolga_book, pagani, kontogiorgis, bruin, dec_jac, elezi_sh, s-sh, deg3, sh_04, sh_03, issac, sh_02, sh_05, fam-gen2, w-pts1, coding-2, beshaj-3 }.

Let $K$ be a genus 3 hyperelliptic field over the ground field $k$. Then $K$ has exactly one genus 0 subfield of degree 2, call it $k(X)$.  It is the fixed field of the \textbf{hyperelliptic involution} $\z_0$ in $\Aut (K)$. Thus, $\z_0$ is central in $\Aut (K)$, where  $\Aut (K)$ denotes the group $\Aut (K/k)$. It induces a subgroup of $\Aut (k(X))$ which is naturally isomorphic to $\bAut (K):= \Aut (K)/\<\z_0\>$. The latter is called the \textbf{reduced automorphism group} of $K$.

 An \textbf{elliptic involution}  of $G = \Aut (K)$  is  an  involution which fixes an elliptic subfield. 
An involution of $\bG=\bAut (K)$ is called \textbf{elliptic} if it is the image of an elliptic involution of $G$.

If $\z_1$ is an elliptic involution in $G$ then $\z_2:=\z_0\, \z_1$  is another involution (not necessarily elliptic). So the   non-hyperelliptic involutions come naturally in (unordered) pairs $\z_1$, $\z_2$. 
These pairs correspond  bijectively to the Klein 4-groups in $G$. 
%Indeed by Remark~\ref{rem0} below,  each Klein 4-group in $G$ contains $\z_0$.
%The latter also correspond to pairs $(E, C)$ of degree 2  subfields of $K$   with   $E\cap k(X)=C\cap k(X)$.

\begin{defn}
We will consider pairs $(K, \e)$  with $K$ a genus 3 hyperelliptic  field and $\e$ an elliptic involution in $\bG$. Two such pairs $(K,\e)$ and  $(K', \e')$ are called isomorphic if there is a $k$-isomorphism  $\a: K \to K'$ with $\e' = \a \e \a^{-1}$.
\end{defn}

Let  $\e$ be an elliptic involution in $\bG$. We can choose the generator $X$ of $\mbox{Fix}(\z_0)$      such that $\e(X)=-X$.  Then $K=k(X,Y)$ where $X, Y$ satisfy equation
\[ Y^2 = (X^2-\a_1^2) (X^2-\a_2^2) (X^2-\a_3^2) (X^2-\a_4^2)  \]
for some $\a_i \in k$, $i=1, \dots, 4$.  Denote by
\begin{equation}
\begin{split}
s_1 = & - \left( \a_1^2 + \a_2^2 + \a_3^2 + \a_4^2 \right) \\
s_2 = & \,  (\a_1 \a_2)^2 + (\a_1 \a_3)^2 + (\a_1 \a_4)^2 + (\a_2 \a_3)^2 + (\a_2 \a_4)^2 + (\a_3 \a_4)^2 \\
s_3 = & -   ( \a_1\,\a_2\,\a_3 )^2 - (\a_4\,\a_1\,\a_2)^2 - (\a_4\,\a_3\,\a_1)^2 - (\a_4\,\a_3\,\a_2)^2 \\
s_4 = & - \left( \a_1 \a_2 \a_3 \a_4 \right)^2\\
\end{split}
\end{equation}
Then,  we have 
\[ Y^2=X^8+s_1 X^6+ s_2 X^4+ s_3 X^2+s_4  \]
with $s_1, s_2, s_3, s_4 \in k$, $s_4\ne 0$. Furthermore,  $E=k(X^2,Y)$  and $ C=k(X^2, YX)$ are the two  subfields corresponding to $\e$ of genus 1 and 2 respectively.   

\par Preserving the condition $\e(X)=-X$ we  can further modify $X$ such that  $s_4=1$. Then, we have the following:

\begin{thm}\label{thm_1}
Let $K$ be a genus 3 hyperelliptic field and $F$ an elliptic subfield of degree 2.  

i) Then,  $K=k(X, Y)$ such that 
\begin{equation}\label{first_eq}
 Y^2 = X^8 + a X^6 + b X^4 +c X^2 +1 
\end{equation} 
for $a, b, c \in k$ such that   $\D \neq 0$, where 
\begin{Small}
\begin{equation}\label{disc_0}
\begin{split}
\D  = &  256(-256+80b^2ac-18ba^3c-18ac^3b-b^2a^2c^2+6a^2c^2-144c^2b+4a^3c^3\\
                   & +4b^3a^2 +4b^3c^2+128b^2-144ba^2+192ac-16b^4+27c^4+27a^4)^2 \\
\end{split}
\end{equation}
\end{Small}

ii) $F = k(U, V) $ where \[ U=X^2,  \quad \textit{ and } \quad V=Y \] and 
\begin{equation}\label{ell_curve}
 V^2= U^4 + a U^3 + b U^2 +c U +1
 \end{equation}

iii) There is a genus 2 subfield $L= k(x, y) $ where  \[ x=X^2,  \quad \textit{ and } \quad y=YX \] and 
\begin{equation}\label{gen_2}
 y^2 = x( x^4+ a x^3 + b x^2 + c  x +1) 
 \end{equation}
\end{thm}

 \proof
 The proof follows from the above remarks.  To show that the genus 2 subfield is generated by $X^2, YX$ it is enough to show that they are fixed by $\z_2$.  
 In cases ii) and iii) we are again assuming that Eq.~\eqref{disc_0} holds.
 \qed
 
%\newpage

\begin{center}
\begin{figure}[htpb]
\[
\xymatrix{                           
                                                                      & K=k(X, Y) \ar@{->}[d]^{\, \, \,2}  \ar@{->}[ld]^{\, \, \,2}  \ar@{->}[rd]^{\, \, \,2}           &                  \\
F=k(X^2, Y)   \ar@{->}[d]^{\, \, \,2}                  & k(X)  \ar@{->}[ld]^{\, \, \,2}  \ar@{->}[rd]^{\, \, \,2}                                                     & L=k(X^2, YX)   \ar@{->}[d]^{\, \, \,2}   \\
k(X^2, Y^2)                                                     &                                                                                                                                   & k( X^2, (YX)^2)       \\
}
\]
\caption{The lattice of a genus 3 hyperelliptic field with an extra involution}
\end{figure}
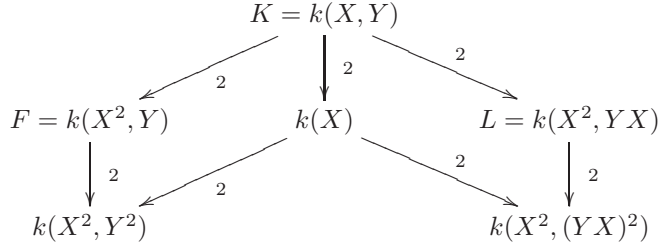
\end{center}

%\newpage
These conditions  determine $X$ up to coordinate change by the group $\< \tau_1, \tau_2\>$ where  
\[ \tau_1: X\to \zz_8 X, \quad \textit{and } \quad   \tau_2: X\to \frac 1 X,\] and  $\zz_8$ is a primitive 8-th root of unity in  $k$. 
Hence,  
\[\tau_1 : \, (a, b, c) \to (\zz_8^6 a, \zz_8^4 b, \zz^2 c ),\]  
and 
\[ \tau_2 : \, (a, b,  c) \to (c, b, a) .\]
Then, $| \tau_1 | =4$ and $|\tau_2 | =2$.  The group generated by $\tau_1$ and $\tau_2$ is the dihedral group of order 8.   Invariants of this action are
\begin{equation}\label{eq2}
\quad \s_2 = \,a\,c,  \quad \s_3=(a^2+c^2)\,b,  \, \quad \s_4=a^4+c^4, 
\end{equation}
since 
\[ 
\begin{split}
& \tau_1 ( a c ) =  \cdot \zz_8^6 a \cdot \zz_8^2 c =  ac \\
& \tau_1 \left(   (a^2 + c^2 ) b \right) = \left( \zz_8^4 a^2 + \zz_8^4 c^2   \right) \cdot ( \zz_8^4 b) = (a^2+c^2) b \\
& \tau_1 (a^4+c^4) = (\zz_8^6 a)^4 + (\zz_8^2 c)^4 = a^4 + c^4 \\
\end{split}
\] 
Since the above transformations are automorphisms of the projective line $\mathbb P^1 (k)$ then the $SL_2(k)$ invariants must be expressed in terms of $\s_2, \s_3$, and $\s_4$.

If  $\s_4 + 2 \s_2^2 =0 $  then this implies that  the curve has automorphism group $\Z_2 \times \Z_4$, continue to Theorem~\ref{thm_last} for details.  From now on we assume that $ \s_4 + 2 \s_2^2 \neq 0$.

In these parameters, the discriminant of the octavic polynomial on the right hand side of
Eq.~\eqref{first_eq} equals 
\begin{small}
\begin{equation}\label{disc}
\begin{split}
\Delta =  &  16\,{s_{{2}}}^{7}+24\,{s_{{2}}}^{6}-72\,s_{{3}}{s_{{2}}}^{5}+16\,s_{{4
}}{s_{{2}}}^{5}+768\,{s_{{2}}}^{5}-1024\,{s_{{2}}}^{4}-2\,{s_{{3}}}^{2
}{s_{{2}}}^{4}+132\,s_{{4}}{s_{{2}}}^{4}\\
& -576\,s_{{3}}{s_{{2}}}^{4}+768
\,s_{{4}}{s_{{2}}}^{3}-72\,s_{{4}}s_{{3}}{s_{{2}}}^{3}+160\,{s_{{3}}}^
{2}{s_{{2}}}^{3}+4\,{s_{{4}}}^{2}{s_{{2}}}^{3}-576\,s_{{4}}s_{{3}}{s_{
{2}}}^{2}+8\,{s_{{3}}}^{3}{s_{{2}}}^{2}\\
& -1024\,s_{{4}}{s_{{2}}}^{2}-{s_
{{3}}}^{2}s_{{4}}{s_{{2}}}^{2}+256\,{s_{{3}}}^{2}{s_{{2}}}^{2}+114\,{s
_{{4}}}^{2}{s_{{2}}}^{2}+192\,{s_{{4}}}^{2}s_{{2}}-18\,{s_{{4}}}^{2}s_
{{3}}s_{{2}}+80\,{s_{{3}}}^{2}s_{{4}}s_{{2}}\\
& -256\,{s_{{4}}}^{2}+27\,{s
_{{4}}}^{3}+128\,{s_{{3}}}^{2}s_{{4}}-16\,{s_{{3}}}^{4}-144\,{s_{{4}}}
^{2}s_{{3}}+4\,{s_{{3}}}^{3}s_{{4}}
\end{split}
\end{equation}

\end{small}

From now forward we will assume that $\D(\s_4, \s_3, \s_2 ) \neq 0$ since in this case the corresponding triple $(\s_2, \s_3, \s_4)$ does not correspond to a genus 3 curve. 
The map $(a, b, c) \mapsto (\s_2, \s_3, \s_4)$ is a branched Galois covering with group $D_4$  of the set 
\[ \{  (\s_2, \s_3, \s_4)\in k^3 : \Delta(\s_2, \s_3, \s_4)\neq 0\}\]
by the corresponding open subset of $(a, b, c)$-space. In any case, it is true that if $a, b, c$ and $a', b', c'$ have the same $\s_2, \s_3, \s_4$-invariants then they are conjugate under $\< \tau_1, \tau_2\>$.

\begin{lemma}\label{lemma1}  For $(a, b, c) \in k^3$ with $\Delta\neq 0$, equation \eqref{first_eq} defines a  genus 3 hyperelliptic  field  $K_{a, b, c}=k(X,Y)$. 
Its reduced automorphism group contains the non-hyperelliptic involution $\e_{a, b, c}: X \mapsto -X$. Two such  pairs $(K_{a, b, c}, \e_{a, b, c})$ and  $(K_{a', b', c'}, \e_{a', b', c''})$ are isomorphic if and only if 
\[ \s_4=\s_4',  \quad \s_3=\s_3', \quad  and \quad \s_2=\s_2', \] 
where $\s_4, \s_3, \s_2$ and $\s_4',\s_3', \s_2'$ are associated with $a, b, c$ and $a', b', c'$, respectively, by   \eqref{eq2}).
\end{lemma}

\proof An isomorphism $\a$ between these two pairs yields $K=k(X,Y)=k(X', Y')$ with $ k(X)=k(X')$ such that $X,Y$ satisfy \eqref{first_eq} and  $X',Y'$ satisfy the corresponding equation with $a, b, c$
replaced by $a', b', c'$. Further, $\e_{a, b, c}(X')=-X'$. Thus $X'$ is conjugate to $X$ under $\< \tau_1, \tau_2\>$ by the above remarks.
This proves the condition is necessary. It is clearly sufficient.

\qed

The following theorem determines relations among $\s_2, \s_3, \s_4$ for each group $G$ such that $V_4 \emb G$.

\begin{thm}\label{thm_last}
Let $\X$ be a  curve in   $\S= \M_3^b \cap \H_3$. Then, the following hold

%i) $\Aut (\X) \iso V_4$

i) $\Aut (\X) \iso \Z_2^3$ if and only if   $\s_4 - 2 \s_2^2 =0 $

ii) $\Aut (\X) \iso \Z_2 \times D_8$ if and only if     $\s_2=\s_4=0$

iii) $\Aut (\X) \iso \Z_2 \times \Z_4$ if and only if   $\s_4+ 2 \s_2^2=0$ and $\s_3 =0$.

iv) $\Aut (\X) \iso D_{12}$ if and only if  
\begin{equation}\label{d_u_12}
\begin{split}
\s_3 & =  {\frac {1}{75}}\, \left( 9\, \s_2 -224 \right)  \left( \s_2-196 \right) \\
\s_4 & = -{\frac {9}{125}}\,{\s_2}^{3}+{\frac {1962}{125}}\,{\s_2}^{2}-{\frac {840448}{1125}}\,\s_2+{\frac {9834496}{1125}}    \\
\end{split}
\end{equation} 

\end{thm}

\noindent See \cite{b-th} for the proof.

\begin{rem}
Any genus 3 hyperelliptic curve is determined uniquely (up to isomorphism) by the tuple $(t_1, \dots , t_6)$ as defined in \cite{hyp_mod_3}.  Such invariants can be expressed in terms of $\s_2, \s_3, \s_4$ as  in \cite{b-th}.  The correctness of the result can be checked by substituting the expressions in the equation of the moduli space $\H_3$, which is computed explicitly in \cite{hyp_mod_3}.
\end{rem}

%************************************************
%\newpage
\section{Decomposition of Jacobians}
%%%%%%%%%%%%%%%%%%%%%%%%%%%%%%%%%%%%%%%%

In this section we give a brief description of the decomposition of Jacobians of genus 3 curves when such decomposition is induced by the automorphisms of the curve. 

Let $\X$ be a genus $g$  algebraic curve with automorphism group $G:=\Aut (\X)$. Let $H \leq G$ such that $H
= H_1 \cup \dots \cup H_t$ where the subgroups $H_i \leq H$ satisfy $H_i \cap H_j = \{ 1\}$ for all $i\neq
j$.  Then,
\[ \Jac^{t-1} (\X ) \times \J^{|H|} (\X / H)\,  \iso \, \J^{| H_1 |} (\X / H_1) \times \cdots \J^{| H_t | } (\X / H_t)\]
The group $H$ satisfying these conditions is called a group with partition. Elementary abelian $p$-groups,
the projective linear groups $PSL_2 (q)$, Frobenius groups, dihedral groups are all groups with partition.

Let $H_1, \dots , H_t \leq G$ be subgroups with $H_i \cdot H_j = H_j \cdot H_i$ for all $i, j \leq t$, and
let $g_{ij}$ denote the genus of the quotient curve $\X/(H_i\cdot H_j)$. Then, for $n_1, \dots , n_t \in \Z$
the conditions
\[ \sum n_i n_j g_{ij} =0, \quad \sum_{j=1}^{t} n_j g_{ij}=0, \]
imply the isogeny relation
\[ \prod_{n_i > 0} \J^{n_i} (\X / H_i) \iso \prod_{n_j < 0} \J^{|n_j|} (\X / H_j)\]
In particular, if $g_{ij}=0$ for $2\leq i < j \leq t$ and if
\[ g = g_{\X/ H_2} + \dots + g_{\X / H_t}\] then
\[ \J (\X) \iso \J (\X /H_2) \times \cdots \times \J (\X / H_t)\]
The reader can check \cite{Ac1, Ac2, KR} for the proof of the above statements.

%******************************************************
\subsection{Non-hyperelliptic curves}
We will use the above facts to decompose the Jacobians of genus 3 non-hyperelliptic curves. $\X$ denotes a
genus 3 non-hyperelliptic curve unless otherwise stated and $\X_2$ denotes a genus 2 curve.

%%%%%%%%%%%%%%%%%%%%%%%%%%%%%%%%%%%%%%%%%%%%%%%%%
\subsubsection{The   group $C_2$} Then the curve $\X $ has an elliptic involution $\s\in \Aut (\X)$. Hence,
there is a Galois covering $\pi \colon \X \to \X/ \< \sigma\>=:E$. We can assume that this covering is maximal.
The induced map $\pi^\ast : E \to \Jac (\X)$ is injective.  Then, the kernel projection $\J (\X) \to E$ is
a dimension 2 abelian variety. Hence, there is a genus 2 curve $\X_2 $ such that
\[ \J (\X_2) \iso E \times \J (\X_2)\]

%**************************
\subsubsection{The Klein 4-group}
Next, we focus on the automorphism groups $G$ such that $V_4 \embd G$.   In this case,   there are three elliptic involutions in $V_4$,
namely $\sigma, \t, \sigma \t$. Obviously they form a partition. Hence, the Jacobian of $\X$ is the product
\[ \J^2 (\X) \iso  E_1^2  \times   E_2^2  \times  E_3^2 \]
of three elliptic curves. By applying the Poincare duality we get
\[ \J (\X) \iso  E_1  \times   E_2  \times  E_3 \]
%

%*********************************
\subsubsection{The dihedral group $D_8$} In this case, we have 5 involutions in $G$ in 3 conjugacy classes.
No conjugacy class has three involutions. Hence, we can pick three involutions such that two of them are
conjugate to each other in $G$ and all three of them generate $V_4$. Hence, $\J (\X) \iso E_1^2 \times E_2$,
for some elliptic curves $E_1, E_2$.

\subsubsection{The symmetric group $S_4$} The Jacobian of such curves splits into a product of elliptic
curves since $V_4 \embd S_4$. Below we give a direct proof of this.

We know that there are 9 involutions in $S_4$, six of which are transpositions. The other three are product
of two  2-cycles and we denote them by $\sigma_1, \sigma_2, \sigma_3$. Let $H_1, H_2, H_3$ denote the subgroups generated
by  $\sigma_1, \sigma_2, \sigma_3$. They generate $V_4$ and are all isomorphic in $G$. Hence, $ \J (\X) \iso E^3,$ for
some elliptic curve $E$.

%**********************************
\subsubsection{The symmetric group $S_3$} We know from above that the Jacobian is a direct product of three
elliptic curves.  Here we will show that two of those elliptic curves are isomorphic. Let $H_1, H_2, H_3$ be
the subgroups generated by transpositions and $H_4$ the subgroup of order 3. Then
\[ \J^3 (\X) \iso E_1^2 \times E_2^2 \times E_3^2 \times \J^3 (\Y)\]
for three elliptic curves $E_1, E_2, E_3$ fixed by involutions and a curve $\Y$ fixed by the element of
order 3. Simply by counting the dimensions we have $\Y$ to be another elliptic curve $E_4$. Since all the
transpositions of $S_3$ are in the same conjugacy class then $E_1, E_2, E_3$ are isomorphic. Then by
applying the Poincare duality we have that
\[ \J (X) \iso E^2   \times E^\prime\]
Summarizing, we have the following:

\begin{thm} Let $\X$ be a genus 3 curve and $G$ its automorphism group. Then,

\noindent a) If $\X$ is hyperelliptic, then the following hold:

i) If $G $ is isomorphic to $V_4$ or $ C_2 \times C_4, $ then $\J (X)$ is isogenous to the product
of and elliptic curve $E$ and the Jacobian of a genus 2 curve $\X_2$,    \[ \J (\X) \iso E \times \J (\X_2).\] 

 ii) If $G $ is isomorphic to $C_2^3$   then $\J (X)$ is isogenous to the product of  three elliptic
curves, \[ \J (\X) \iso E_1 \times E_2 \times E_3.\] 

iii) If $G $ is isomorphic to $D_{12},   C_2 \times S_4 $ or any of the groups of order 24 or 32, then $\J
(X)$ is isogenous to the product of three elliptic curve such that two of them are isomorphic  \[ \J (\X) \iso E_1^2 \times E_2.  \]

\noindent b) If $\X$ is non-hyperelliptic then the following hold: 

i) If $G $ is isomorphic to $C_2$ then $\J (X)$ is isogenous to the product of an elliptic curve and the
Jacobian of some genus 2 curve $\X_2$
\[ \J (\X) \iso E \times \J (\X_2)\]
ii) If $G $ is isomorphic to $V_4$ then $\J (X)$ is isogenous to the product of three elliptic curves
\[ \J (\X) \iso E_1 \times E_2 \times E_3\]
iii) If $G $ is isomorphic to $S_3, D_8$ or has order 16 or 48 then $\J (X)$ is isogenous to the product of
three elliptic curves such that two of them are isomorphic to each other
\[ \J (\X) \iso E_1^2 \times E_2 \]
iv) If $G $ is isomorphic to $S_4, L_3(2)$ or $C_2^3 \sem S_3$ then $\J (X)$ is isogenous to the product of
three elliptic curves such that all three of them are isomorphic to each other
\[ \J (\X) \iso E^3 .\]
\end{thm}

\proof The proof of the hyperelliptic case is similar and we skip the details. The reader interested in
details can check   \cite{P}.

Part b): When $G$ is isomorphic to $C_2, V_4,  D_8, S_4, S_3$ the result follows from the remarks above. The
rest of the theorem is an immediate consequence of the list of groups as in the Table 1 of \cite{g-sh}. If $|G|=16, 48$ then $D_8 \embd   G$.  Then, from
the remarks at the beginning of this section the results follows. If $G$ is isomorphic to $L_3(2)$ or $C_4^2
\sem S_3$ then $S_4 \embd G$. Hence the Jacobian splits as in the case of $S_4$. This completes the proof.
\endproof

It is possible that given the equation of $\X$ one can determine the equations of the elliptic or genus 2  components in all cases of the theorem.

The reader interested in decomposition of Jacobians of curves with large automorphism group can check \cite{lange}, where an algorithm is provided. Indeed, the method suggested in \cite{lange} works for all curves $C$ which have an element $\tau \in \Aut (C)$ such that $C/\< \tau \>$ has genus 0.  This includes all  superelliptic curves.  The automorphism groups of such curves (including positive characteristic) were determined in \cite{Sa2} and their equations in \cite{Sa1}.   
%A complete decomposition of Jacobians of superelliptic curves is intended in \cite{dec_jac}. 

The above theorem gives the splitting of the Jacobian based on automorphisms.  
Next we will focus on the $(2, 4, 4)$ splitting for hyperelliptic curves.  We will explicitly determine the elliptic components for a given genus 3 curve $\X$.  

%*******************************************************************
%\newpage
 \section{Subcovers of genus 2}
 
 Next we study in detail the complement $C$ of $E$ in $\Jac (\X)$.  From the above theorem, $C$ has equation as in  Eq.~\eqref{gen_2}.  Its absolute invariants $i_1, i_2, i_3$, as defined in \cite{sh-v},  can be expressed in terms of the dihedral invariants $\s_4, \s_3, \s_2$ as follows:
\begin{small}
\begin{equation}\label{i1_i2_i3}
\begin{split}
i_1  = &  \,  {144} \frac M {D^2} \,  \left( 2\,{\s_2}^4 +2\,{\s_2}^3 - 6\,{\s_2}^2 \s_3+{\s_2}^2 \s_4 - 40\,{\s_2}^2 + \s_2 \s_4 + 9\,{\s_3}^2 - 3\,\s_3 s_4 -20\,\s_4    \right)  \\ \\
i_2  = & \,      432 \, \frac {M^2} {D^3} \, \left( 8\,{\s_2}^5  +12\,{\s_2}^4    -36\,{\s_2}^3  \s_3+4\, {\s_2}^3  \s_4  +1248\,{\s_2}^3-558\,{\s_2}^2\s_3+114\,{\s_2}^2\s_4 \right. \\
& \left. +81\ ,\s_2 {\s_3}^2 -18\,\s_2 \s_3 \s_4 -2240\,{\s_2}^2+624\,\s_2\s_4 +216\,{\s_3}^2-279\,\s_3\s_4 +54\,{\s_4}^2-1120\,\s_4 \right)     \\ \\
i_3  =  & \, {\frac {243} {16}    \,  \frac {M^3} {D^5}}\,  \left( 16\,{\s_2}^{7}+24\,{\s_2}^{6}-72\,{\s_2}^{5}\s_3+16\,{\s_2}^{5}\s_4-2\,{\s_2}^4  {\s_3}^2+768\,{\s_2}^{5}-576\,{\s_2}^{4}  \s_3   \right. \\
& \left.  +132\,{\s_2}^4  \s_4 +160\,{\s_2}^3{\s_3}^2-72\,{\s_2}^3\s_3\s_4  +4\,{\s_2}^3{\s_4}^2+8\,{\s_2}^2{\s_3}^3-{\s_2}^2{\s_3}^2\s_4   -1024\,{\s_2}^{4}  \right. \\
& \left. +768\,{\s_2}^3 \s_4 +256\,{\s_2}^2 {\s_3}^2-576\,{\s_2}^2\s_3\s_4  +114\,{\s_2}^2{\s_4}^2    +80\,\s_2{\s_3}^2\s_4    -18\,\s_2\s_3{\s_4}^2-16\,{s_3}^{4} \right. \\
& \left. +4\,{\s_3}^3\s_4  -1024\,{\s_2}^2\s_4  +192\,\s_2{\s_4}^2+128\,{\s_3}^2\s_4  -144\,\s_3{\s_4}^2+27\,{\s_4}^3-256\,{\s_4}^2
  \right)    \\ \\
\end{split}
\end{equation}
\end{small}
where $M= \s_4 + 2\s_2^2$ and $D= 16\,{s_{{2}}}^{3}-40\,{s_{{2}}}^{2}+8\,s_{{2}}s_{{4}}-3\,{s_{{3}}}^{2}-20\,s_{{4}}$.  For the rest of the paper we assume that  $D=J_2 \neq 0$.

We are interested in the case when the Jacobian of $C$ splits.  Splitting of such Jacobians of genus 2 has been studied in our previous work at \cite{sh-v, deg3}.  We only consider the case when $\Jac (C)$ is $(2, 2)$ decomposable.  The locus $\L_2$ of such genus two curves is computed in \cite{sh-v} in terms of the invariants $i_1, i_2, i_3$.  Substituting the expressions in Eq.~\eqref{i1_i2_i3} in the equation of $\L_2$ from \cite{sh-v} we have the following:
\begin{equation}\label{eq_T}
\left( 2\,{\s_{{2}}}^{2}-{\s_{{4}}} \right)  \cdot \left( 2\,{\s_{{2}}}^{2} + {\s_{{4}}} \right) \cdot F_1 (\s_2, \s_3,  \s_4)  = 0
\end{equation}
where $F_1 (\s_2, \s_3,  \s_4)$ is an irreducible polynomial of degree 13, 8, 6 in given in $\s_2, \s_3,  \s_4$ respectively.   The expression of $F_1 (\s_2, \s_3,  \s_4)$ is given in  the Appendix. 

Let the locus of genus 3 hyperelliptic curves whose Jacobian is $(2, 4, 4)$-split  be denoted by $\mathcal T$.    There are three components of $\mathcal T$ which we denote them by $\mathcal T_i$, $i=1, 2, 3$ as seen by Eq.~\eqref{eq_T}.

Two of these components are well known and the correspond to the cases when $V_4$ is embedded in the reduced  automorphism group of $\X$.  These cases correspond to the singular locus of $\S_2$ and are precisely the locus $\det \left( \Jac (\phi) \right)=0$, see \cite{sh_03}.  This happens for all genus $g\geq 2$ as noted by Geyer \cite{Ge}, Shaska/V\"olklein \cite{sh-v}, and shown in \cite{sh_03}.

\begin{lemma}
Let $\X$ be a genus 3 curve with $(2, 2, 4)$-split Jacobian.  Then, one of the following occurs

i) $ \Z_2^3 \embd \Aut (\X)$

ii) $\Z_2 \times \Z_4  \embd \Aut (\X)$ 

iii) $\X$ is in the locus $\mathcal T_3$ 
\end{lemma}

\proof  The proof is an immediate consequence of Theorem 2 and Eq.~\eqref{eq_T}.   

\qed

The third component $\mathcal T_3$  is more interesting to us.  
%It doesn't seem to have any group theoretic reason for this component to be there in the first place.   
It is the moduli space of pairs of degree 4 non-Galois covers $\psi_i: \X_3 \to E_i$, $i=1, 2$. 
 
One of the main goals of this paper is to determine explicitly the family $\mathcal T_3$ of genus 3 curves and relations among its elliptic subcovers.
We have the maps 
\begin{equation}
\begin{aligned}
  \mathcal T_3  &  \overset{\psi}{\longrightarrow} \L_2    \overset{\psi_0}{\longrightarrow}  k^2  \\
  (\s_2, \s_3, \s_4 )  &\to  (u, v)  \to (j_1, j_2) \\
\end{aligned}
\end{equation}
where $\s_2, \s_3, \s_4$ satisfy  $F_1 (\s_2, \s_3, \s_4)=0$  and $u, v$ are given explicitly by Eq.~\eqref{i1_i2_i3} and Thm.~ (3) in \cite{sh-v}. 
%The degree $\deg \, \psi_0 =2$ and $\deg \, \psi =39 $.

The expressions of $u$ and $v$ are computed explicitly in terms of $\s_2, \s_3, \s_4$ by substituting the expressions of Eq.~\eqref{i1_i2_i3} expressions for $u$ and $v$ as rational functions of $i_1, i_2 i_3$ as computed in   Shaska/V\"olklein \cite{sh-v}. As rational functions $u$ and $v$ have degrees 35 and 70 respectively (in terms of $\s_2, \s_3, \s_4$). The $j$-invariants $j_1$ and $j_2$ will be determined in the next section. 

Since $\mathcal T_3$ is a subvariety of $\H_3$ it would be desirable to express its equation in terms of a coordinate in $\H_3$.  One can use the absolute invariants of the genus 3 hyperelliptic curves $t_1, \dots , t_6$ as defined in \cite{hyp_mod_3} and the expressions of $\s_4, \s_3, \s_2$ in terms of these invariants as computed in \cite{b-th}. 

\begin{rem}  $\mathcal T_3$ is a 2-dimensional subvariety of $\H_3$ determined by the equations 
\begin{equation}\label{eq_T_3} 
\left\{ 
\begin{split}
F_1 (\s_2, \s_3, \s_4 ) =0 \\
t_i - T_i (\s_2, \s_3, \s_4) ,  \quad i = 1, \dots , 6 \\
\end{split}
\right.
\end{equation}
where $T_i$ is the function $t_i$ evaluated for the triple $(\s_2, \s_3, \s_4)$. 
\end{rem}

The equations of $\mathcal T_3$ can be explicitly determined in terms of $t_1, \dots , t_6$ by eliminating $\s_2, \s_3, \s_2$ from the above equations. Normally, when we talk about $\mathcal T_3$ we will think of it given in terms of $t_1, \dots , t_6$.

\begin{exa}[Automorphism group $\Z_2^3$]
Consider the genus 3 curves $\X$ with $\Aut (\X) \iso \Z_2^3$.  Then, $\s_4 = 2 s_2^2$ and 
\[
\begin{split}
u  = & \frac 1 P \, \left(  -9\,{s_{{3}}}^{2}+120\,s_{{2}}s_{{3}}-400\,{s_{{2}}}^{2}+16\,{s_{{2}}}^{3} \right) \\
v  = &  - \frac {2}     { P^2   }  \, 
( 432\,s_{{2}}{s_{{3}}}^{3} -27\,{s_{{3}}}^{4} -1440\,{s_{{2}}}^{2}{s_{{3}}}^{2}-6400\,{s_{{2}}}^{3}s_{{3}}     +32000\,{s_{{2}}}^{4}+288\,{s_{{2}}}^{3}{s_{{3}}}^{2} \\
& -5376\,{s_{{2}}}^{4}s_{{3}}+23040\,{s_{{2}}}^{5}+256\,{s_{{2}}}^{6}    )   \\
\end{split}
\]
where $P =  -s_3^2 -8\,s_2 s_3    -16\, s_2^2 +16\, s_2^3$. 
\end{exa}

%\newpage

 For the rest of this section we will see if we can invert the map $\psi$.  Given an ordered pair $(u, v) \in \L_2$ there is a unique genus two curve corresponding to this pair.  Indeed, an explicit equation can be found in terms of $u$ and $v$ from the following.

\begin{prop}\label{thm:model-racional-uv}
Let $(u,v)\in k^2$ such that
  $$(u^2-4v+18u-27)(v^2-4u^3)(4v-u^2+110u-1125)\neq 0.$$
Then, the curve of genus $2$ defined over $k$ given by
\begin{equation}\label{rat_v4}
y^2=a_0 x^6+a_1 x^5+a_2 x^4+a_3 x^3+t a_2 x^2+t^2 a_1 x+t^3 a_0,
\end{equation}
corresponds to the moduli point  $(u,v) \in \L_2 \emb \M_2$,  where one of the following holds:
  
i)  If   $u\neq 0$,  then
  \begin{align*}
    t&=v^2-4u^3,\\
    a_0&=v^2+u^2v-2u^3,\\
    a_1&=2(u^2+3v)(v^2-4u^3),\\
    a_2&=(15v^2-u^2v-30u^3)(v^2-4u^3),\\
    a_3&=4(5v-u^2)(v^2-4u^3)^2,
  \end{align*}

ii)  If $u=0$,  then
\[    t=1,\, a_0=1+2v,\,     a_1=2(3-4v),\,     a_2=15+14v,\,    a_3=4(5-4v)\]
\end{prop}

A proof of this is provided in \cite{min_eq} together with other computational aspects of genus 2 curves. 

Hence, corresponding to the pair $(u, v)$ there is a unique genus 2 curve $C_{u, v}$.  The following Lemma addresses the rest of our question. 
 
\begin{lemma}\label{lem_gen_2}
i) Any genus two curve $C$ defined over an algebraically  closed field $k$ can be written as 
\begin{equation}\label{g3-g2} y^2 =  x( x^4+a x^3+b x^2+c x +1) \end{equation}
for some $a, b, c \in k$ such that 
\begin{small}
\[ 
\begin{split}
\Delta(a, b, c) = & \, \, 256-192 c a-6 c^2 a^2-4 c^3 a^3-27 c^4-27 a^4-80 b^2 c a-128 b^2\\
                  &+18 b c a^3+18 c^3 a b+144 b a^2+144 c^2 b+c^2 b^2 a^2-4 b^3 a^2-4 c^2 b^3+16 b^4 \neq  0 
\end{split}
\]
\end{small}

ii) Let $C$ be a genus 2 curve with equation as in Eq.~\eqref{g3-g2}.  Then, there exists a genus 3 curve $\X$ with equation 
\[ Y^2 = X^8 + aX^6 + bX^4 + CX^2 +1 \]
and a degree 2 map $f : \X \to C$ such that 
\[ x= X^2 \quad \textit{and } \quad y = YX\]
\end{lemma}

\proof
i) Let $C$ be a genus 2 curve defined over $k$.  Then, the equation of $C$ is given by \[ y^2 = \Pi_{i=1}^6 (x-\a_i), \] 
where $\a_i$ are all distinct for all $i=1 \dots 6$.  Since $k$ is algebraically closed, then we can pick a change of transformation in $\P$ such that $\a_1 \to 0$ and  $\a_2 \to \infty$.  We can also pick a coordinate such that $\a_3 \cdots a_6 =1$.  Then, the curve $C$ as equation as claimed. The condition that $\Delta (a, b, c) \neq 0$ simply assures that not two roots of the sextic coalesce. 

ii)  This genus 3 curve is a covering of $C$ from Thm.~\ref{thm_1}. 

\qed

Hence, the curve $C_{u, v}$ can be written as in Eq.~\eqref{g3-g2}.  This would mean that we can explicitly compute $\s_4, \s_3, \s_2$ in terms of $u$ and $v$.  Finding a general formula for $(\s_2, \s_3, \s_4)$ in terms of $(u, v)$ is computationally difficult. Under some additional restrictions this ca be done, as we will see in the next section.

%******************************************************************
%\newpage
\section{Elliptic subfields}

In this section we will determine the elliptic subcovers of the genus 3 curves $\X \in \mathcal T_3$.  We will describe how this can be explicitly done, but will skip displaying the computations here.

A point $\p =(t_1, \dots , t_6) \in \mathcal T_3$ satisfies equations Eq.~\eqref{eq_T_3}. Our goal is to determine the $j$-invariants of $E, E_1, E_2$ in terms of $t_1, \dots t_6$.  The $j$-invariant of $E$ is 
\begin{equation}\label{eq_j}
\begin{split}
j = & 256\,  \frac { \left( -{\s_{{3}}}^{2}-12\,\s_{{4}}-24\,{\s_{{2}}}^{2}+3\,\s_{{2}}\s_{{4}}+6\,{\s_{{2}}}^{3} \right) ^{3}}
{ \left( \s_{{4}}+2\,{\s_{{2}}}^{2} \right)  }
( 4\,{\s_{{3}}}^{3}\s_{{4}}     -256\,{\s_{{4}}}^{2}-1024\,{\s_{{2}}}^{4} \\
& -72\,{\s_{{2}}}^{5}\s_{{3}}-144\,{\s_{{4}}}^{2}\s_{{3}} +8\,{\s_{{2}}}^{2}{\s_{{3}}}^{3}-576\,{\s_{{2}}}^{4}\s_{{3}}+114\,{\s_{{2}}}^{2}{\s_{{4}}}^{2}+132\,{\s_{{2}}}^{4}\s_{{4}}-2\,{\s_{{2}}}^{4}{\s_{{3}}}^{2}\\
& +160\,{\s_{{2}}}^{3}{\s_{{3}}}^{2}+4\,{\s_{{2}}}^{3}{\s_{{4}}}^{2}+16\,{\s_{{2}}}^{5}\s_{{4}}-16\,{\s_{{3}}}^{4}+192\,\s_{{2}}{\s_{{4}}}^{2}+128\,{\s_{{3}}}^{2}\s_{{4}}-1024\,{\s_{{2}}}^{2}\s_{{4}}\\
& +768\,{\s_{{2}}}^{3}\s_{{4}}+24\,{\s_{{2}}}^{6}+16\,{\s_{{2}}}^{7}+27\,{\s_{{4}}}^{3}+768\,{\s_{{2}}}^{5}-18\,\s_{{2}}{\s_{{4}}}^{2}\s_{{3}}-72\,{\s_{{2}}}^{3}\s_{{4}}\s_{{3}}\\
& -576\,{\s_{{2}}}^{2}\s_{{3}}\s_{{4}}-{\s_{{2}}}^{2}{\s_{{3}}}^{2}\s_{{4}}+80\,\s_{{2}}{\s_{{3}}}^{2}\s_{{4}}+256\,{\s_{{2}}}^{2}{\s_{{3}}}^{2} ) \\
\end{split}
\end{equation}  

We denote the degree 2 elliptic subcovers of $C$ by $E_1$ and $E_2$ and their $j$-invariants by $j_1$ and $j_2$.  These $j$ invariants are the roots of the quadratic
\begin{equation}\label{eq_j1_j2}
 j^2+   256 \,  \frac { 2 u^3-54 u^2+9 u v-v^2+27 v } { u^2+18 u-4 v-27 } \, j + 65536 \, \frac {u^2+9 u-3 v} {(u^2+18 u-4 v-27)^2}, 
 \end{equation}
see Eq.~(4) in \cite{sh-v}. 

\begin{center}
\begin{figure}[htpb]
\[
\xymatrix{                           
                                                                      & \X   \ar@{->}[ld]_{\, \, \,2}  \ar@{->}[rd]^{\, \, \,2}           &                  \\
E               &                                   & C   \ar@{->}[ld]_{\, \, \,2} \ar@{->}[rd]^{\, \, \,2}  &  \\
                                                  &                                                                                                                     E_1              &   &    E_2\\
}
\]
\caption{A genus 3 curve $\X$ with its three elliptic subcovers}
\end{figure}
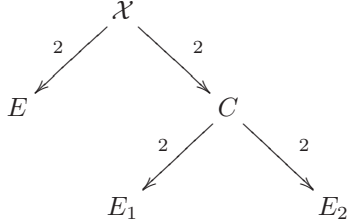
\vspace{-.5cm}
\end{center}
Since these $j$-invariants are determined explicitly in terms of $u$ and $v$, then via the map $\psi: \mathcal T_3 \to \L_2$ we express such coefficients in terms of $\s_2, \s_3, \s_4$.  
Moreover, the maps 
\begin{equation}
\begin{split}
 & \mathcal T_3 \to k^3\setminus \{\D(\s_2, \s_3, \s_4) = 0\} \to k^2\setminus\{\D_{u, v} = 0\} \to k^2 \\
& (t_1, \dots , t_6) \to (\s_4, \s_3, \s_2 ) \overset{\psi}{\longrightarrow} (u, v) \to (j_1, j_2) \\
\end{split}
\end{equation} 
are all  explicitly determined.

\begin{rem}
Given a tuple $(t_1, \dots , t_6)$ which satisfies the equations of $\mathcal T_3$ does not necessarily represent a genus 3 hyperelliptic curve.  We have to check that such tuple satisfies the equation of the genus 3 hyperelliptic moduli, as described in \cite{hyp_mod_3}. 
\end{rem}

\begin{exa}
Let be given a 6-tuple \[(t_1, \dots , t_6) = \left( {\frac {8767521}{6224272}}, {\frac {152464}{5329}},   {\frac {8116}{3431}},  -{\frac {3343532}{695617}}, -{\frac {91532}{148117}},  -{\frac {50448727768}{28398241}} \right) \]
which satisfies the Eq.~(17) in \cite{hyp_mod_3}.  Then, this tuple correspond to a genus 3 hyperelliptic curve, more precisely 
the  curve $\X$ with equation \[ Y^2= X^8 + X^6 + X^4 + X^2+1\]
Then, the corresponding invariants are  \[ \s_2 =1, \quad  \s_3=2, \quad s_4=2\]
The genus 2 subcover has invariants 
\[ i_1 = - \frac {48} 5, \quad i_2= \frac {432} 5, \quad i_3 = \frac 1 {400} \]
and the corresponding dihedral invariants 
$u$ and $v$   are 
\[ u = 9, \qquad v = - \frac {754} 5 \]
The $j$-invariants of the three elliptic subcover are 
\[ j= 2048, \quad j_1= {\frac {32768}{5}}+ \frac 2 5\,\sqrt {268435081} , \quad j_2 = {\frac {32768}{5}} - \frac 2 5\,\sqrt {268435081} \]
\qed
\end{exa}

Next we will study a very special subvariety of $\mathcal T_3$, such as the locus of curves in $\mathcal T_3$ when  $E_1$ is isomorphic  to $E_2$. 

%*************************************
%\newpage
\subsection{Isomorphic elliptic subfields}

The two elliptic curves $E_1$ and $E_2$  are isomorphic when their $j$-invariants are equal, which happens when the discriminant of the  quadratic in Eq.~\eqref{eq_j1_j2} is zero.  From Remark~(1) in \cite{sh-v} this occurs if and only if 
\[ (v^2-4u^3) (v-9u+27) =0\]
The first condition is equivalent to $D_8 \embd Aut (C)$. The later condition gives 
\[ u = 9 - \frac \lambda {256}, \quad \textit{and } \quad v = 9 \left( 6 - \frac \lambda {256} \right), \]
where $\lambda:=j_1=j_2$.    Both of these loci can be explicitly computed given enough computing power.

Substituting $u$ and $v$ in terms of $\s_2, \s_3, \s_4$ in the equation \[ v-9u+27 =0, \]
we get an equation of degree 68, 42, and 29 in $\s_2, \s_3, \s_4$ respectively.  We denote it by
\begin{equation}\label{f_2}
 F_2 (\s_2, \s_3, \s_4) =0 
\end{equation} 
and do not display it here because of its size.      This equation and the Eq.~\eqref{eq_T}  define the locus $\mathfrak T$ in terms of $\s_4, \s_3, \s_2$.

%The rational expressions of $\s_4, \s_3, \s_2$ in terms of the absolute invariants $t_1, \dots , t_6$ are computed in \cite{b-th}.  Hence, we have a way of computing $\mathfrak T$ in terms of $t_1, \dots , t_6$. 

Let $\mathfrak T$ denote the algebraic variety defined by by the equations 
\begin{equation}
\left\{
\begin{split}
& F_1(\s_2, \s_3, \s_4) =0\\
& F_2(\s_4, \s_3, \s_2 ) =0\\
\end{split}
\right.
\end{equation}

\begin{lemma}
The algebraic variety $\mathfrak T$ is a 1-dimensional subvariety, it has 5 genus 0 components as in Eq.~\eqref{T_comp}. 
Every point $(\s_2, \s_3, \s_4) \in \mathfrak T$ correspond to a genus 3 hyperelliptic curve with $(2, 4, 4)$-split Jacobian such that the degree 4 elliptic subcovers are isomorphic to each other. 
\end{lemma}

\proof From the equations above we can eliminate $\s_3$ via resultants and get the following. 
In this case we get
\[ (2\s_2^2 - \s_4 )^{16}  (\s_4 + 2 \s_2^2)^{172} \,  g_1^{12} \, g_2^{12} \, g_3^{10} \, g_4^8 \, g_5 =0\]
where 
\begin{Small}
\begin{equation}\label{T_comp}
\begin{split}
g_1  = & \, s_{{4}}+2\,{s_{{2}}}^{2}-100\,s_{{2}}+625  \\
g_2  = & \, -27\,s_{{4}}+{s_{{2}}}^{3}+6\,{s_{{2}}}^{2}+768\,s_{{2}}-4096 \\
g_3  = & \, -16777216+5242880\,s_{{2}}-450560\,{s_{{2}}}^{2}+7680\,{s_{{2}}}^{3}-340\,{s_{{2}}}^{4}+8\,{s_{{2}}}^{5}-102400\,s_{{4}} \\
& +16640\,s_{{2}}s_{{4}}-220\,{s_{{2}}}^{2}s_{{4}}+4\,s_{{4}}{s_{{2}}}^{3}-125\,{s_{{4}}}^{2}   \\
g_4  = & \, 3515625-937500\,s_{{2}}+62500\,{s_{{2}}}^{2}+64\,{s_{{2}}}^{4}+15000\,s_{{4}}-2000\,s_{{2}}s_{{4}}  \\
g_5  = & \, 35153041-16173862\,s_{{4}}+2926323\,{s_{{4}}}^{2}+131244344\,s_{{2}}-11788512\,{s_{{2}}}^{4}+133698380\,{s_{{2}}}^{2}\\
& +41050\,{s_{{4}}}^{3}+1792\,{s_{{2}}}^{7}+22764\,{s_{{4}}}^{2}{s_{{2}}}^{3}-21768\,{s_{{2}}}^{5}s_{{4}}+136952\,{s_{{2}}}^{6}-134208\,{s_{{2}}}^{5}\\
& +16762008\,{s_{{2}}}^{3}-57934\,{s_{{2}}}^{2}{s_{{4}}}^{2}-111744\,s_{{4}}{s_{{2}}}^{4}-1040380\,s_{{2}}{s_{{4}}}^{2}+1987608\,s_{{4}}{s_{{2}}}^{3}+2559786\,{s_{{2}}}^{2}s_{{4}}\\
& -27622208\,s_{{2}}s_{{4}}-16\,{s_{{2}}}^{7}s_{{4}}-316\,{s_{{2}}}^{6}s_{{4}}+16\,{s_{{2}}}^{5}{s_{{4}}}^{2}-4\,{s_{{2}}}^{3}{s_{{4}}}^{3}-530\,{s_{{2}}}^{2}{s_{{4}}}^{3}-6180\,s_{{2}}{s_{{4}}}^{3}\\
& +685\,{s_{{2}}}^{4}{s_{{4}}}^{2}+132\,{s_{{2}}}^{8}+125\,{s_{{4}}}^{4}  \\
\end{split}
\end{equation}
\end{Small}

We have  $ (2 \s_4-\s_2^2)  (2 \s_4+\s_2^2) \neq 0$, as noted before. All other components  are  genus zero curves.  

\qed

\begin{rem}
The equation of $\mathfrak T$ can be expressed in the absolute invariants $t_1, \dots , t_6$ by eliminating $\s_4, \s_3, \s_2$ from expressions in Eq.~\eqref{i1_i2_i3}.  Such expressions are large and we do not display them here. 
\end{rem}

Let be given a parameterization of $\mathfrak T$.  Then we have the following maps
\begin{equation}
\begin{split}
k  \to \mathfrak T & \to L_2 \to k \\
t   \to \left( \s_4(t), \s_3(t), \s_2(t) \right) & \to \left( u(t), v(t) \right) \to j(t) \\
\end{split}
\end{equation}

This map gives us the possibility to construct a family of curves defined over $\Q$ such that all their subcovers, namely $C$, $E$, $E_1$, and $E_2$ are also defined over $\Q$.  For example, for $t \in \Q$ we have the corresponding $\s_4, \s_3, \s_3 \in \Q$.  Hence, there is a genus 3 curve $\X$ defined over $\Q$.  The invariants $u, v$ are rational functions of $\s_4, \s_3, \s_2$ and therefore of $t$.  Thus, $u, v \in \Q$.  Form Prop.~\ref{thm:model-racional-uv} there a genus 2 curve $C$ such that $C$ is defined over $\Q$.  Moreover, the $j$-invariants for all elliptic subcovers are rational functions in $\s_4, \s_3, \s_2$ and therefore in $t$.  Hence, $E$, $E_1$, $E_2$ are also defined over $\Q$.

\begin{thm}
Let $\X$ be a  curve in   $\mathfrak T$ and $\s_2, \s_3, \s_4$ its corresponding dihedral invariants.  Then 
\[ \J (\X) \iso E \times E^\prime \times E^\prime \]
where $E$ and $E^\prime$ are elliptic curves with  $j$-invariants $j(E)$ as in Eq.~\eqref{eq_j} and $j(E^\prime )$ as 
\[
\begin{split}
j^\prime & = -128\,{\frac {2\,{u}^{3}-54\,{u}^{2}+9\,uv-{v}^{2}+27\,v}{{u}^{2}+18\, u-4\,v-27}},\\
\end{split}
\]
where $u$ and $v$ are given as rational functions of $i_1, i_2, i_3$ as in \cite{sh-v}. 
Moreover, there is only a finite number of genus 3 curves $\X$ such that  $E \iso E^\prime$.  
%This occurs if 
%
%\[ G(\s_2, \s_3, \s_4) =0.  \]
\end{thm}

\proof
The equation of $j(E)$ was computed in Eq.~\eqref{eq_j}.  Since the other two elliptic subcovers have the same $j$-invariants then this invariant is given by the double root of the quadratic in Eq.~\eqref{eq_j1_j2}. Thus, 
\[ j^\prime =   -128\,{\frac {2\,{u}^{3}-54\,{u}^{2}+9\,uv-{v}^{2}+27\,v}{{u}^{2}+18\, u-4\,v-27}}   \]
Substituting the values for $u$ and $v$ we get the expression as claimed. 

We have  $E \iso E^\prime$ if and only if $j = j^\prime$.  This gives a third equation $G(\s_2, \s_3, \s_4) =0$ as claimed in the theorem. By Bezut's theorem, the number of solutions of the system of equations $F_i (\s_2, \s_3, \s_4) =0$, for  $i=1, 2$ and $G(\s_2, \s_3, \s_4)=0$  is finite. 

\qed

The family above could be significant in number theory in constructing genus 3 curves with many rational points. Similar techniques have been used for genus 2 in \cite{sh_02} and by various other authors.

%***************************************
\appendix
\section*{Appendix}
\renewcommand{\thesection}{A}

The polynomial $F_1 (\s_2, \s_3, \s_4)$ from Eq.~\eqref{eq_T}. 

\begin{small}
\begin{equation*} 
\begin{split}
F_1 &(\s_2, \s_3 \s_4)  = 1024\,{\s_2}^{13}-256\,{\s_2}^{12}\s_3+1536\,{\s_2}^{12}-8448\,{\s_2}^{11}\s_3+2560\,{\s_2}^{11}\s_4+1664\,{\s_2}^{10}{\s_3}^2 \\
& -512\,{\s_2}^{10}\s_3\s_4+64\,{\s_2}^{9}{\s_3}^3+33600\,{\s_2}^{11}-48720\,{\s_2}^{10}\s_3+10752\,{\s_2}^{10}\s_4+37696\,{\s_2}^{9}{\s_3}^2\\
& -20928\,{\s_2}^{9}\s_3\s_4+2560\,{\s_2}^{9}{\s_4}^2-4448\,{\s_2}^{8}{\s_3}^3+3008\,{\s_2}^{8}{\s_3}^2\s_4-384\,{\s_2}^{8}\s_3{\s_4}^2-432\,{\s_2}^{7}{\s_3}^4 \\
& +96\,{\s_2}^{7}{\s_3}^3\s_4-4\,{\s_2}^{6}{\s_3}^{5}+20000\,{\s_2}^{10}-140000\,{\s_2}^{9}\s_3+84000\,{\s_2}^{9}\s_4+172600\,{\s_2}^{8}{\s_3}^2\\
& -139200\,{\s_2}^{8}\s_3\s_4+21120\,{\s_2}^{8}{\s_4}^2-71680\,{\s_2}^{7}{\s_3}^3+78368\,{\s_2}^{7}{\s_3}^2\s_4-20736\,{\s_2}^{7}\s_3{\s_4}^2 \\
& +1280\,{\s_2}^{7}{\s_4}^3+2288\,{\s_2}^{6}{\s_{{3}}}^4-6200\,{\s_2}^{6}{\s_3}^3\s_4+2016\,{\s_2}^{6}{\s_3}^2{\s_4}^2-128\,{\s_2}^{6}\s_3{\s_4}^3+1256\,{\s_2}^{5}{\s_3}^{5}\\
& -568\,{\s_2}^{5}{\s_3}^4\s_4+48\,{\s_2}^{5}{\s_3}^3{\s_4}^2+26\,{\s_2}^4{\s_3}^{6}-4\,{\s_2}^4{\s_3}^{5}\s_4+50000\,{\s_2}^{8}\s_4+80000\,{\s_2}^{7}{\s_3}^2\\
& -280000\,{\s_2}^{7}\s_3\s_4+84000\,{\s_2}^{7}{\s_4}^2-140000\,{\s_2}^{6}{\s_3}^3+345500\,{\s_2}^{6}{\s_3}^2\s_4-156600\,{\s_2}^{6}\s_3{\s_4}^2\\
& +19200\,{\s_2}^{6}{\s_4}^3+30000\,{\s_2}^{5}{\s_3}^4-107520\,{\s_2}^{5}{\s_3}^3\s_4+61008\,{\s_2}^{5}{\s_3}^2{\s_4}^2-10272\,{\s_2}^{5}\s_3{\s_4}^{3}\\
& +320\,{\s_2}^{5}{\s_4}^4+8060\,{\s_2}^4{\s_3}^{5}-204\,{\s_2}^4{\s_3}^4\s_4-2628\,{\s_2}^4{\s_3}^3{\s_4}^2+592\,{\s_2}^4{\s_3}^2{\s_4}^3-16\,{\s_2}^4\s_3{\s_4}^4 \\
& -1464\,{\s_2}^3{\s_3}^{6}+1256\,{\s_2}^3{\s_3}^{5}\s_4-244\,{\s_2}^3{\s_3}^4{\s_4}^2+8\,{\s_2}^3{\s_3}^3{\s_4}^3-72\,{\s_2}^2{\s_{3}}^{7}+21\,{\s_2}^2{\s_3}^{6}\s_4 \\
& -{\s_2}^2{\s_3}^{5}{\s_4}^2+50000\,{\s_2}^{6}{\s_4}^2+120000\,{\s_2}^{5}{\s_3}^2\s_4-210000\,{\s_2}^{5}\s_3{\s_4}^2+42000\,{\s_2}^{5}{\s_4}^3+40000\,{\s_2}^4{\s_3}^4 \\
& -210000\,{\s_2}^4{\s_3}^3\s_4+259350\,{\s_2}^4{\s_3}^2{\s_4}^2-87000\,{\s_2}^4\s_3{\s_4}^3+9120\,{\s_2}^4{\s_4}^4+30000\,{\s_2}^3{\s_3}^4\s_4 \\
& -53760\,{\s_2}^3{\s_3}^3{\s_4}^2+21080\,{\s_2}^3{\s_3}^2{\s_4}^3-2544\,{\s_2}^3\s_3{\s_4}^4+32\,{\s_2}^3{\s_4}^{5}-3600\,{\s_2}^2{\s_3}^{6}+8060\,{\s_2}^2{\s_{3}}^{5}\s_4\\
& -1920\,{\s_2}^2{\s_3}^4{\s_4}^2-202\,{\s_2}^2{\s_3}^3{\s_4}^3+64\,{\s_2}^2{\s_3}^2{\s_4}^4-732\,\s_2{\s_3}^{6}\s_4+314\,\s_2{\s_3}^{5}{\s_4}^2-34\,\s_2{\s_3}^4{\s_4}^3\\
& +81\,{\s_3}^{8}-36\,{\s_3}^{7}\s_4+4\,{\s_3}^{6}{\s_4}^2+25000\,{\s_2}^4{\s_4}^3+60000\,{\s_2}^3{\s_3}^2{\s_4}^2-70000\,{\s_2}^3\s_3{\s_4}^3+10500\,{\s_2}^3{\s_4}^4\\
& +40000\,{\s_2}^2{\s_3}^4\s_4-105000\,{\s_2}^2{\s_3}^3{\s_4}^2+86525\,{\s_2}^2{\s_3}^2{\s_4}^3-23925\,{\s_2}^2\s_3{\s_4}^4+2208\,{\s_2}^2{\s_4}^{5}\\
& +7500\,\s_2{\s_3}^4{\s_4}^2-8960\,\s_2{\s_3}^3{\s_4}^3+2728\,\s_2{\s_3}^2{\s_4}^4-252\,\s_2\s_{{3}}{\s_4}^{5}-1800\,{\s_3}^{6}\s_4+2015\,{\s_3}^{5}{\s_4}^2\\
& -623\,{\s_3}^4{\s_4}^3+59\,{\s_3}^3{\s_4}^4+6250\,{\s_2}^2{\s_4}^4+10000\,\s_2{\s_3}^2{\s_4}^3-8750\,\s_2\s_3{\s_4}^4+1050\,\s_2{\s_4}^{5}+10000\,{\s_3}^4{\s_4}^2 \\
& -17500\,{\s_3}^3{\s_4}^3+10825\,{\s_3}^2{\s_4}^4-2610\,\s_3{\s_4}^{5}+216\,{\s_4}^{6}+625\,{\s_4}^{5} \\
\end{split}
\end{equation*}
\end{small}

%*********************************************************************************************************
\bibliographystyle{amsplain}

\end{document}